\documentclass[12pt]{article}
\usepackage{amssymb, url, amsthm, amsmath}

\usepackage{textcomp}
\def\boxit#1{%
{\hbox{\lower3pt\hbox{\vrule\vbox{\hrule\kern2pt%
\hbox{\kern2pt$#1$\kern2pt}\kern2pt\hrule}\vrule}}}}
\bigskip
\bigskip
\def\be{\begin{equation}}
\def\ee{\end{equation}}

\def\R{{\sf I\kern-.2em R}}
\def\N{{\sf I\kern-.2em N}}
\def\C{\kern.1em{\raise.47ex\hbox{$\scriptscriptstyle
$}}\kern-.40em{\sf C}}
\def\Z{{\sf Z\kern-.32em Z}}

\def\hat{\widehat}

\def\hat{\widehat}
\def\be{\begin{equation}}
\def\ee{\end{equation}}

\newtheorem{theorem}{\noindent Theorem}

\newtheorem{definition}{\noindent Definition}
\newtheorem{corollary}{\noindent Corollary}

\newtheorem{statement}{\noindent Proposition}

\urldef\vershik\url{vershik@pdmi.ras.ru}

\author {A.~M.~Vershik\thanks{%
St.~Petersburg Department of Steklov Institute of Mathematics.
E-mail: \vershik. Partially supported by the RFBR grant
 05-01-00899 and  CRDF grant RUM1-2622-ST-04.%
}}

\date{16.02.07}
\title{Krein duality, positive 2-algebras, and dilation of comultiplications}

\begin{document}
 \maketitle
\rightline{\it To the centenary of Mark Krein}

\begin{abstract}
The Krein--Tannaka duality for compact groups
was a generalization the Pontryagin--Van Kampen duality for locally
compact abelian groups and a remote predecessor of the theory of
tensor categories. It is less known that it
found applications in algebraic combinatorics (``Krein algebras'').
Later, this duality was substantially extended: in \cite{V},
the notion of {\it involutive algebras in positive vector duality}
was introduced. In this paper, we
reformulate the notions of this theory
using the language
of bialgebras (and Hopf algebras) and introduce the class of
involutive bialgebras and positive 2-algebras.
The main goal of the paper is to give a precise statement of
a new problem, which we consider as one of the main
problems in this field,
concerning the existence of
dilations (embeddings) of positive 2-algebras into involutive bialgebras,
or, in other words, the problem of describing subobjects
in involutive bialgebras. We define two types of subobjects in the category of
bialgebras, strict and nonstrict ones, and consider the corresponding
embeddings (dilations)
of positive 2-algebras into bialgebras. The difference between the two
types of dilations is illustrated by the example of bicommutative
positive 2-algebras (commutative hypergroups). The most interesting
instance of our problem concerns dilations of the Hecke algebra $H_n(q)$.
It seems that in this case strict dilations may exist only
for $q=p^k$ (with $p$ a prime); it is not known whether nonstrict dilations
exist for other $q$. We also show that the class of finite-dimensional involutive
semisimple bialgebras coincides with the class of semigroup algebras of
finite inverse semigroups.
\end{abstract}

 \section{Introduction}

The goal of this paper is to formulate the problem of lifting
an involutive algebra with positive multiplication and comultiplication
(we call such an object a positive 2-algebra) into an involutive
bialgebra, for example, into the group bialgebra of a finite group
or an inverse semigroup. This problem arose in connection with
some problems of the theory of association schemes, and also
in connection with some combinatorial problems of representation theory.
But the dilation problem itself is akin to some
questions of operator theory (dilations of semigroups of operators) and
the theory of dynamical and algebraic systems.
In short, we call it
the problem of dilating comultiplications or
derandomizing a random multiplication.

Many ideas in this field spring from M.~G.~Krein's works on
duality and positivity; that is why we begin (Sec.~2)
with a brief historical survey on positive multiplications and
duality, which shows the complicated way in which M.~G.~Krein's ideas
on duality and positivity of algebras
have finally united nowadays. This survey
is not complete and does not touch upon many close questions
(see \cite{Va}).

The most impressive concrete instance of our main problem is given in
Sec.~3; it concerns the Hecke algebra $H_n(q)$
and its embedding into $GL(n,\mathbb C)$. We consider this problem from different
points of view.

In Sec.~4 we define the category of {\it positive involutive
2-algebras}, i.e., spaces with positive (with respect to the given
involution and coinvolution) multiplication and comultiplication,
which will be called positive 2-algebras.
The category of such 2-algebras includes the category of bialgebras
\cite{Ca} and Hopf algebras: in this case, the multiplication and the comultiplication
are not homomorphisms, but only positive maps.
Functorially, it is equivalent to the category of {\it $*$-algebras
in positive duality} introduced by the author \cite{V} in 1971,
but in this paper we use the language of bialgebras and multiplications-comultiplications instead of the language of algebras in duality.
Let us emphasize that the idea of duality itself is in a sense
a continuation of M.~G.~Krein's ideas. The problem of lifting
positive operations to multiplicative ones is traditional
for functional analysis and its applications (see \cite{FN}). Here we formulate it
in a purely algebraic context. In Sec.~5 we define two kinds of subobjects in
involutive bialgebras, and then (Sec.~6) give a precise statement
of the general problem, defining what it means to embed
a positive 2-algebra into a bialgebra. For want of space, we omit
examples and calculations for some (bicommutative) cases, see \cite{Ve}.
The general dilation problem is far
from being solved.

\section{M.~G.~Krein's work on duality and positivity in algebras}

In the 40s, M.~G.~Krein worked in representation theory and published
several papers on this subject. The papers  \cite{K1,K2}
are devoted to invariant Hermitian-positive kernels on homogeneous spaces
and are close to his favourite circle of problems concerning positive definiteness
and its relations to operator theory.\footnote{``I respond to the words
`a positive definite
function' like a war-horse to a trumpet call,''
M.G.\ said, in the elevated style typical for him, in his talk at a conference on
functional analysis held in Odessa in 1958.}
As far as I know, he did not continue the work started in
these papers; nevertheless, they became
rather widely known. In particular, in what follows we will need
one of the results from
\cite{K2}, a theorem on the equivalence of algebraic and scalar
positivity. In the pre-war paper
\cite{K3} and the more complete paper \cite{K4} on the same subject,
which seem to be in no way connected to the previous ones
(though the paper \cite{K4} cites the main result of \cite{K2}),
the Pontryagin--Van Kampen duality for locally
compact abelian groups, very popular at that time,
was extended to compact nonabelian groups
(somewhat earlier this had been done, though in less generality,
by T.~Tannaka; M.G.\ learned about this from D.~A.~Raikov after the
paper  \cite{K4} had been published).

While in the abelian case, the duality does not lead out of the category
of abelian groups: the ``Fourier functor'' associates with a given group
the group of its characters and generalizes the classical
Fourier transform theory from $\Bbb Z$ and $\Bbb R$
to all locally compact abelian groups, in the nonabelian case,
the dual object is no longer a group, but an algebra.
What is now called a Krein block algebra \cite{K1}
is a commutative algebra partitioned into an at most countable sum of simple
finite-dimensional algebras over  $\Bbb C$ with a distinguished basis.
This basis is the basis of matrix units in all classes of equivalent
irreducible representations of the group. More precisely, if $G$ is
a compact group, then one can easily see that its Krein-dual block
algebra is the algebra of all continuous functions
on $G$ endowed with a decomposition into the direct sum of
minimal two-sided ideals and a {\it basis} consisting of the matrix elements
of all unitary irreducible representations. The Krein duality theorem
says that this block algebra determines the group up to an isomorphism:
it is the group of scalar homomorphisms that are multiplicative on the basis.
The proof exploits, in a spectacular way, Gelfand's theory of
normed rings.\footnote{I dare to conjecture that it is the success of
the recently appeared Gelfand's theory of normed rings, and especially
the theory of infinite-dimensional representations, that stimulated
M.~G.~Krein and A.~N.~Kolmogorov to work on these problems (see their
correspondence in \cite{Ko}).} M.G.\ realizes the group algebra of
a compact group as, in his words, a block algebra, i.e., an algebra
of functions on the countable set of classes of irreducible unitary
representations with matrix values; choosing an appropriate basis
in this algebra, we arrive at the duality theorem.

A more symmetric statement of Krein's theorems arises if we get rid of
bases and use the invariant language of the theory of duality of algebras
as vector spaces, see Sec.~5. The first algebra is the group algebra of the group $G$
over some field;
and the dual algebra is the algebra
of continuous functions on the group with values in the same field. Such a formulation,
in the spirit of Bourbaki's duality theory for vector spaces \cite{Bour},
was suggested by the author in \cite{V}, where he introduced the notion
of pairs of $*$-algebras in duality and its geometric analog, the theory
of packets. This theory was elaborated in
\cite{Ke1, Ke2} by S.~V.~Kerov and, more recently, in \cite{VEP,BVEP,E}.

The modern formulation of such a duality, developed gradually
starting from the papers of that time, uses the terminology of bialgebras
and Hopf algebras. It is as follows: in the group algebra there is a
commutative comultiplication, in the commutative algebra of functions
on the group there is a convolution comultiplication, and each of these two
bialgebras determines the other one. This formulation, without any
modifications, works for finite groups; later it was generalized to
universal enveloping algebras of Lie groups. But in the general case
of locally compact groups, it requires serious topological comments,
and M.G.'s paper in fact provided those considerations
from the theory of normed rings that suffice for compact groups.
In continuation of the works of Tannaka and Krein, duality
for general  locally compact groups was considered by
many other authors. The theory of duality of group algebras and algebras of
functions on groups eventually resulted, though not only for internal reasons,
in the theory of quantum groups. Besides, these works provided a basis
for a number of direct generalizations: Tannaka--Krein categories in algebraic
geometry  \cite{De, Sa}, Tannaka-Krein duality for groupoids, the notion
of monoidal categories
(see \cite{Mas, VDa, Jo-Str}), the theory of hypergroups
and multivalued groups,
the theory of generalized shift,
noncommutative
integration (I.~Segal, W.~F.~Stinespring,
G.~I.~Kac, and others). All these topics are considered in several hundred papers
and several dozen books.

However, it should be noted that duality for noncommutative groups, when it is
formulated in abstract form and does not involve some special properties
of the group,
is tautological; moreover, I dare say that it
had no serious analytic consequences comparable to the Fourier theory
for abelian groups. The cause is that the category dual to the category
of group algebras is difficult to describe independently in invariant terms,
even in the finite-dimensional case. In order to construct an interesting
generalization of Fourier theory and harmonic analysis, one needs an additional
structure on the group or the algebra
under consideration, which fixes a certain commutative subgroup
(subalgebra) defined in invariant terms, for example, the Cartan subgroup
(in the theory of semisimple Lie algebras), the Gelfand--Tsetlin algebra
(in the theory of symmetric and similar groups), etc. Then the duality becomes
an analytic tool for studying groups or algebras.

At the same time, the true development of Krein's ideas on duality
took quite another course, in which these ideas were combined with
the positivity theorems from \cite{K1, K2}. Namely,
a real application, which apparently was not quite foreseen by the authors,\footnote{In 1980,
M.G.\ showed me a book on algebraic coding theory and relation schemes sent
him by the author (as far as I remember, P.~Delsarte), where his name
appeared on
every two pages, and said with some embarrassment, ``I have nothing to do with this
and do not understand this at all.'' ``This'' was precisely the Krein
duality, but not for the group algebras of compact groups, as in his
papers, but for algebras generated by association schemes.}
 was discovered much later, on the one hand, in the theory of duality
of algebras (see above), and on the other hand, with significant delay,
in algebraic combinatorics and the theory of association schemes
(see, e.g.,  \cite{Ban}). It is in this field that the notion of  a
block (Krein) algebra became a model for generalizations
and started to be used not only for group algebras. In modern terms,
this was the transition from bialgebras and Hopf algebras to
systems (2-algebras) with a {\it less rigid agreement between
multiplications and comultiplications}. And the crucial argument
for constructing this theory was precisely one
of the positivity theorems from
\cite{K2} rather than duality itself.

In short, the intricate chain of references to these papers by M.~G.~Krein looks
as follows. The positivity condition for a very special case was introduced
by D.~Higman in 1964. Then L.~Scott published a brief note in {\it Notices}
\cite{Sc1} saying that his colleague C.~Dunkl,
a specialist in harmonic analysis,
drew his attention to the fact that Krein's paper  \cite{K1}
contains a general positivity condition for homogeneous spaces, so that
Higman's condition, as Scott joked, is new only for the example
considered by Higman. This was followed by papers of
D.~Higman \cite{Hi}, L.~Scott \cite{Sc2}, and many others on the
algebraic theory of association schemes (see  \cite{Ban}),
after which the condition of positivity of the Krein constants became
generally accepted and widely used. It is interesting that all
these papers, as well as subsequent papers on this subject,
are strongly influenced
by the paper \cite{K4} on block algebras, and it is the notion
of a block algebra that they generalize to a nongroup situation
(which was not considered
by M.G.), but in fact they use not so much the results of
\cite{K4} on constructing the block algebra of a compact group, as
the simple positivity condition for functions
on homogeneous spaces,
which plays a crucial role in all
these considerations, from the paper \cite{K2}. Note that the preliminary
paper \cite{K3} published in 1941 was noticed by S.~Bochner
\cite{B} already
during the war; apparently, that is the reason why this result of M.G.\
became known in the West rather early, at least to specialists in functional analysis.
It is also interesting that these papers by M.~G.~Krein were rediscovered
and exploited in the studies on algebraic
combinatorics, with almost 30 years delay, by mathematicians working in the United States and Japan
rather than Ukraine or Moscow, where there were many specialists
working on
exactly the same problems and where
M.~G.~Krein's works must have been more known. However, when preparing the above-mentioned paper
\cite{V}, in which the so-called positive duality of $*$-algebras was introduced
and the positivity of multiplication and comultiplication was postulated,
I also was guided rather by the duality theorems of Krein and Tannaka
and their followers
than by M.G.'s papers \cite{K1, K2, K4}, which I already knew.
The paper \cite{V}
on positive pairings of $*$-algebras was conceived as an attempt to
extend the Krein duality from groups to homogeneous spaces
of double cosets, but then it turned out that the suggested weakening of the multiplicativity
condition covers much more general objects, not
necessarily of group origin. But in this paper, we
relate the main topic
to another field, in which M.G.\ also worked, namely, to the theory
of dilations. This link is apparently new. In functional analysis,
the theory of dilations is represented
by works of Foias and Nagy, Krein's school,
W.~F.~Stinespring
(the theorem on dilation of completely positive maps), M.~A.~Naimark
(lifting of positive operator measures to projection measures or
spectral measures), etc.; now it appears also in problems of
purely algebraic nature.

\section{An instance of the main problem}

\subsection{Dilation of the classical Hecke algebra}

The complex Hecke algebra $H_n(q)$, $q\in (0,\infty)$, $n\geq 2$, $n\in{\Bbb N}$,
is the algebra with generators $\tau_i$,  $i=1,\dots, n-1$, satisfying the relations
   $$\tau_i \tau_j =\tau_i \tau_j, \qquad |i-j|>1,\; i,j=1, \dots, n-1,$$
   $$(\tau_i+1)(\tau_i-q)=0,\qquad i=1,\dots, n-1,$$
   $$\tau_i \tau_{i+1}\tau_i=\tau_{i+1}\tau_i \tau_{i+1},  \qquad i=1, \dots,
   n-2. $$
It is well known that this algebra is isomorphic to the group algebra
${\Bbb C}({\frak S}_n)$ of the symmetric group ${\frak S}_n$, and the elements
   $\tau_i$, $i=1,\dots, n-1$, are deformations of the classical Coxeter
   transpositions $\sigma_i$. A linear basis of the Hecke algebra is indexed
   by the elements of the symmetric group:
   ${\frak S}_n \ni g \rightarrow \tau_g$, where the Coxeter transposition
   $\sigma_i=(i,i+1)\in {\frak S}_n $ is associated to the generator
   $\tau_{\sigma_i}\equiv\tau_i$, $i=1, \dots ,n-1$,
   and an arbitrary permutation $g$ with {\it reduced} decomposition
   into a product of Coxeter transpositions  $g=\prod_k \sigma_{i_k}$
   is associated to the element $\tau_g=\prod_k \tau_{i_k}$.
   Consider the coordinatewise cocommutative comultiplication in this algebra,
   which is diagonal in the above basis
      $\{\tau_g\}$:
   $$\Delta(\tau_g) =\tau_g \bigotimes \tau_g,\qquad g\in {\frak S}_n.$$

For $q \ne 1$, the linear basis $\{\tau_g\}$, $g\in {\frak S}_n $,
does not form a group (the product of two elements of this basis,
    $\tau_{h_1}\cdot \tau_{h_2}=\sum_g c_{h_1,h_2}^g \tau_s$,
is only a linear combination of basis elements), so that
{\it this comultiplication is no longer a homomorphism of the algebra}
   $H_n(q)$ to the algebra $H_n(q)\bigotimes H_n(q)$. Hence the Hecke
   algebra is not a bialgebra, and {\it a fortiori} not a Hopf algebra,
   which is the case for the group algebras
  ${\Bbb C}({\frak S}_n)$. At the same time, if $q=p^m$ is a power of
  a prime number $p$, then a classical result says that the Hecke
  algebra   $H_n(q)$ is a subalgebra of the group bialgebra ${\Bbb C}(GL_n(F_q))$
  of the group $GL_n(F_q)$ of all invertible matrices
  over the finite field $F_q$, namely, the subalgebra
  of double cosets of a Borel subgroup (the group of upper triangular matrices);
  in this case, the comultiplication in
  $H_n(q)$ introduced above is a projection of the comultiplication
  in the group algebra. Thus the Hecke algebra
  $H_n(q)$ for $q=p^m$ can be embedded into the group algebra of the group
   $GL_n(F_q)$.

   However, for $q>1$ {\it the comultiplication in the Hecke algebra,
   being no longer a homomorphism, still remains a positive map}
   in the sense of the cone of positive elements (with respect to
   the natural involution); indeed, as can easily be seen, the coefficients
   $c_{g,h}^s$, the structure constants of the multiplication, are nonnegative.
   It is obvious that the positivity of the comultiplication and
   multiplication is a necessary condition for the embeddability
   into a bialgebra, since homomorphisms are always positive maps and
   the positivity is preserved under positive projections (but, of course,
   not every positive map is multiplicative). Using the language of operator theory,
   we may say that for $q=p^m$ there exists a
   {\it dilation} (lifting) of the Hecke algebra: the original positive but
   not multiplicative operations, the multiplication and comultiplication
   in $H_n(q)$, in the ambient group bialgebra
   ${\Bbb C}(GL_n(F_q))$ are restrictions and projections of multiplicative
   operations.\footnote{Examples of using the term ``dilation'' are
   as follows: a unitary dilation of a contraction
   $T$ in a Hilbert space $K$ is a unitary operator $U$ in an extension
   $H\supset K$ such that for the orthogonal projection
   $P:H\to K$ and all positive integers $n$, the equality
   $T^n=PU^n$ holds (Naimark-type or Nagy-type theorems, etc.);
   or another situation, which is closer to the one under consideration:
   a dilation of a positive (Markov) operator $T$ in $L^2$
   is a unitary multiplicative operator $U$ with a similar property:
   $T^n=PU^n$ for some expectation $P$ (see \cite{VPOL}).
   In our case, we consider only a dilation of the comultiplication,
   regarded as an operation from the algebra to its tensor square,
   since the multiplication is directly inherited from the ambient
   bialgebra. Dilations were also considered by M.~G.~Krein
   and Foias and Nagy, see \cite{FN}. In contrast to known theorems on dilations
   of a single operator, we want to find, for example, coherent dilations
   of the family of operators of multiplication by all the basis elements
   of the algebra.}

   Thus a natural question arises, whether the positivity of the
   multiplication and comultiplication is a sufficient condition for
   the embeddability into some bialgebra. More precisely,

\smallskip
   {\it whether a similar embedding also exists for other values
      $q\ne p^m$, i.e., whether the Hecke algebra
   $H_n(q)$ for other values of $q$ can be embedded, with the same multiplication
   and a dilation of the comultiplication, into some bialgebra,
   not necessarily a group one, or into a Hopf algebra, even perhaps
   infinite-dimensional? Indeed, the necessary condition
   for the existence of such an embedding, the positivity
   of the multiplication and comultiplication, is satisfied.}
\smallskip

To state the problem, we need to specify what the term
``embedding'' or ``dilation'' means; to this end, we define two
notions of a subobject in a bialgebra.
The preliminary analysis of these notions
and a partial solution of the problem for two-dimensional positive 2-algebras
is   the main   subject of this paper.

The existence of a strict (see below) embedding of the Hecke algebra $H_n(q)$
into the {\it group algebra of
   a finite group} for general $q$ is unlikely; for instance, for
$n=3$, $q\ne p^m$ embeddings into group bialgebras do not exist (an observation
   due to I.~N.~Ponomarenko). On the contrary, it is rather
   the existence of such a group ($GL_n(q)$) for $q=p^m$ that looks surprising.
And though, as we will see, the class of involutive semisimple bialgebras is not
exhausted by group bialgebras,  and the complete description of
this class given here shows that it coincides with the class of semigroup bialgebras of
finite inverse semigroups, nevertheless this does not increase the
possibilities for the positive answer.

A description of {\it obstacles for the existence of such embeddings} would be
   not less interesting for the analysis of such objects. They appear
   already in dimension~2. An obvious obstacle for the existence of
   {\it finite-dimensional} covering algebras is the irrationality
   of the structure constants of the comultiplication. An obstacle
   for the existence of {\it group dilations} is a violation of certain
   arithmetic relations between the structure constants; alternatively, such an obstacle
   may be related to conditions on the orders of groups and subgroups, etc.
Anyway, the author does not know the
general answer concerning the existence of
dilations of the Hecke algebras into bialgebras.

  \subsection{Probabilistic interpretation:
  derandomization of a random multiplication}

   Our example with the Hecke algebra can be viewed in another way.
   Let us make a change of generators:
      $$\bar{\tau}_g=q^{-l(g)}\tau_g, \qquad g \in {\frak S}_n,$$
   where $l(g)$ is the reduced length of a permutation $g$; in particular,
   ${\bar\tau}_i=\frac{\tau_i}{q}$, ${\bar e}=e$.
   Now consider the formula
 $${\bar \tau}_{h_1}\cdot{\bar\tau}_{h_2}=\sum_g \mu_{h_1,h_2}^g {\bar \tau}_g$$
   for the multiplication of elements of the new basis. It is easy to check
   that in the chosen normalization
   the sum $\sum_s \mu_{h_1,h_2}^g$
   becomes equal to one  (and all these coefficients are nonnegative), so that we may say that the product of any two
   elements of the new basis is a probability measure
   $\mu_{h_1,h_2}$ on the set $\{\tau_g,\, g \in {\frak S}_n\}$. In other words,
   the multiplication is ``probabilistic,'' or ``random,'' in contrast
   to the multiplication in groups or semigroups, which is ``deterministic''
   (in a clear sense).

Given a group or a semigroup, or a more
general system with a binary operation, let us say that a partition of this system
into subsets (blocks) is {\it stable} if the vector subspace of all linear combinations
of blocks is an algebra. Then the above question concerning
a dilation of the multiplication can be formulated as follows. Assume that
we are given a ``random'' multiplication on the set ${\frak S}_n$
(or on another set); thus the vector space of functions on this set
is an algebra. The question is whether the multiplication in this algebra is a
projection (randomization) of a deterministic multiplication?
In more detail, does there exist a group or a semigroup (or even a $G$-space)
and a stable partition of this system
such that for any two blocks of this partition
($s$ and $t$), the distribution of their product, regarded as a measure
on the set of blocks, coincides with the multiplication
rule in the original algebra (the measure $\mu_{s,t}$ in the example above)?
In this interpretation, our question is
{\it how to derandomize a random multiplication or comultiplication}.
Analogous but simpler derandomization problems arise, for instance, in
the general theory of Markov processes, ergodic theory, etc. Such a derandomization
is suggested at the end of the paper.
Another, combinatorial, reformulation of our problem is as follows:
is it possible to represent the multiplication and comultiplication tensors
of a positive 2-algebra (see below) with rational real coefficients
as factors of multiplication and comultiplication tensors with
coefficients $(0,1)$? This aspect of the problem will be considered
separately.

\section{Bialgebras, involutive bialgebras}

\subsection{Bialgebras}

\begin{definition}
A \textbf{2-algebra} is a vector space over
$\Bbb C$ in which there is a structure of an associative algebra with
unit and that of a coassociative coalgebra with counit
over the same field (in general, without any conditions
on the compatibility of these structures). \textbf{In what follows,
unless otherwise stated, we consider only finite-dimensional 2-algebras
that are semisimple algebras and coalgebras (a coalgebra is semisimple
if its dual algebra is semisimple).}
A bialgebra (see {\rm \cite{Ca}}) is a 2-algebra satisfying the following
equivalent conditions: {\rm1)} the comultiplication and the counit are homomorphisms
of the corresponding algebras; {\rm2)} the multiplication and the unit
are homomorphisms of the corresponding coalgebras.\footnote{It would be
natural to use the term ``bialgebra'' for objects that
we have called 2-algebras, and use an epithet to express
one or another kind of compatibility between the multiplication
and the comultiplication; then what is
now called a bialgebra should be called a ``multiplicative bialgebra,''
what is now called a ``positive 2-algebra'' (see below) should be called a
``positive bialgebra,'' etc. But, unfortunately,
it is too late for  such a change of
terminology.}
\end{definition}

Denote the multiplication, regarded as an operator from
$A$ to $A\bigotimes A$, by
$\delta$, the unit by $\epsilon$, the comultiplication by $\Delta: A\rightarrow
A\bigotimes A $, and the counit by $\varepsilon$. The definition
of a bialgebra is symmetric with respect to the pairs
$(\Delta,\varepsilon)$ and
$(\delta,\epsilon)$, which often makes it  superfluous to
prove parallel assertions.
Sometimes we will also write the product in
$A$ in the ordinary way: $\delta(x,y)=x\cdot y$.

It is useful to relate all these and subsequent definitions
with the definitions of pairs
of finite-dimensional algebras in duality. Let $A$ be
a 2-algebra. Recall that in the finite-dimensional case the dual space
to an algebra (coalgebra) is a coalgebra (respectively, an algebra).
Let us define a multiplication in the dual vector space $A'$
as the map
adjoint to the comultiplication $\Delta$, an involution in $A'$ as the
operation adjoint to the coinvolution in $A$, and a unit in $A'$
as the counit in $A$. Then we may speak of the pair of algebras $A$ and $A'$
in duality
with some or other additional properties (multiplicativity, positivity, etc.).
Conversely, an ordered pair of algebras in a nondegenerate duality
generates a 2-algebra.
This observation will be used in what follows.

The group algebra of a finite group with convolution multiplication
and diagonal comultiplication is obviously a cocommutative bialgebra
(and even a Hopf algebra). It is also well known (see \cite{Ca}) that
the semigroup algebra of any finite semigroup (monoid) with unit,
equipped  with the natural operations,
is also a cocommutative bialgebra. In what
follows, we will speak of {\it group} or {\it semigroup}
bialgebras, meaning that they are endowed with
both multiplicative structures.

\subsection{Involutive bialgebras}

Our plan is to introduce an additional structure in bialgebras, which
somehow or other exists in group bialgebras and bialgebras close to them, and then to
weaken some requirements on the structures of bialgebras.
Let us equip a bialgebra $A$ with an involution and a coinvolution;
an involution in
$A$ is a second-order antilinear antiautomorphism of the algebra $A$;
similarly,
a second-order antilinear antiautomorphism of the
coalgebra $A$ is called a coinvolution. Hopf algebras with involution
(but without coinvolution) were considered earlier; see, e.g.,
\cite{Ca}. Denote the involution and the coinvolution in $A$ by
$\sharp$ and $\flat$, respectively. In this notation, the relations between
the involution, coinvolution, multiplication, and comultiplication
look as follows:
 $$\delta (\sharp\otimes\sharp) =\sharp \delta J, \qquad (\flat \otimes \flat)\Delta =J \Delta \flat,$$
where $J$ is the permutation (flip) in the tensor product:
 $ J(x\otimes y)=y\otimes x$.

\begin{definition} A bialgebra equipped with an involution and a coinvolution
is called an \textbf{involutive bialgebra} or a bialgebra with involution
if the multiplication commutes with the coinvolution and the comultiplication
commutes with the involution:
$$\Delta\sharp=(\sharp\otimes\sharp)\Delta, \qquad
\delta(\flat\otimes\flat)=\flat\delta.$$
\end{definition}

Given an algebra with an involution $\sharp$, we can define the notion of
{\it nonnegative elements}, i.e., elements of the form
$x\cdot x^{\sharp}$. They form a convex {\it cone}
$K^{\sharp}$, which will be called the
{\it positive cone with respect to the involution $\sharp$}
(in the finite-dimensional case, this cone will be closed).
{\it The positive
cone $K^{\flat}$ with respect to the coinvolution $\flat$}
is defined as follows: this is the complete preimage under the comultiplication
$\Delta$ of the conic hull of the set of elements of the form
$x\otimes x^{\flat}\in
A\bigotimes A$, $x\in A $; in other words, $K^{\flat}=\{y\in A: \Delta
y=\sum_i a_i\otimes a_i^{\flat}\}$.\footnote{In terms of duality,
$K^{\flat}$ is just the cone dual to the cone of positive elements
in the dual algebra with respect to the involution adjoint
to $\flat$.}

\begin{theorem}
The comultiplication in an involutive bialgebra preserves the positive cone
with respect to the involution, and the multiplication preserves
the positive cone with respect to the coinvolution; in short,
in an involutive bialgebra the multiplication and the comultiplication are
\textbf{positive}.
\end{theorem}

\smallskip\noindent{\bf Remark.}
This assertion contains, in particular, the algebraic part of Krein's
theorem on positivity of multiplication and comultiplication for the
finite-dimensional case.
\smallskip

\begin{proof}
We must prove that the coproduct
of an element nonnegative with respect to the involution
lies in the
positive cone, and that the product of elements nonnegative with respect to the
coinvolution lies in the corresponding cone. By symmetry, it suffices
to prove one of these claims; the second one will follow from the first
one by replacing the comultiplication with the multiplication and the involution
with the coinvolution, and passing to the dual bialgebra. When proving
that the coproduct of positive elements is positive, we may restrict
ourselves to considering elements of the form
$x\cdot x^{\sharp}$, $x\in A$, since they generate the cone of positive elements.
The following calculations use the commutation relations stated above:
 $$\Delta(x\cdot x^{\sharp})=\Delta(x)\Delta (x^{\sharp})
 =\Delta(x)[\Delta(x)]^{\sharp}.$$
We have used the multiplicativity of the comultiplication
and the fact that it
commutes with the involution; similarly, the proof of the second
claim uses the multiplicativity of the multiplication (with respect to
the comultiplication) and the fact that it
commutes with the coinvolution.
\end{proof}

\subsection{Involutive bialgebras and semigroup algebras
of inverse semigroups}

The complex group algebras of finite groups with the natural structures are
involutive bialgebras in the sense defined above: the involution
is generated by taking the inverse element in the group
($g \to g^{-1}$), and the coinvolution is generated by taking the conjugates of the
coefficients. Moreover, an involutive cocommutative Hopf algebra
is the group bialgebra of a finite group. Below we will consider the class
of finite-dimensional involutive semisimple cocommutative bialgebras and
relate it to an important class of semigroups.

Namely, we will show that the class of finite inverse unital semigroups
generates exactly the class of involutive semisimple bialgebras.
Recall (see \cite{KP}) that an inverse semigroup is a semigroup $S$ in which
every element $a \in S$ has a unique inverse element,
i.e., an element
$b \in S$ such that $aba=a$.
The uniqueness of such an element guarantees that the inverse element to $b$ is $a$,
so that in an inverse semigroup there is a natural involution, which associates
with each element $a$ its inverse element denoted by
$a^*=a^{-1}$. The multiplication, comultiplication, and coinvolution
are defined in the same way as in the group case. If the semigroup $S$
has a two-sided identity element, then the functional ``the value at
this identity element''
is a counit (in the coalgebra).

The main example of an inverse semigroup, which is most interesting for us, is
the {\it symmetric inverse semigroup}
${\cal I}_n$
of one-to-one partial maps of an $n$-element set into itself
(i.e., bijections between subsets of this set),
with the
empty map as the zero of the semigroup. An important theorem, which
generalizes Cayley's theorem on groups, says that every finite inverse semigroup
has an isomorphic embedding into a symmetric inverse semigroup.

\medskip\noindent{\bf Example.}
Consider the inverse semigroup ${\cal I}^1_n$
consisting of the $n^2$ matrix units
   $\{e_{i,j}\}$, $i,j=1, \dots, n$, and the zero $0$
with the ordinary matrix multiplication.  The inverse element is defined
as follows:   $e_{i,j}^*=e_{j,i}$, $0^*=0$. The semigroup algebra
${\Bbb C}({\cal I}^1_n)$ is the direct sum $M_n(\Bbb C)\bigoplus \Bbb C$
of the algebra of matrices
and the one-dimensional
two-sided ideal, equipped with the Kronecker
comultiplication (i.e., in the
dual formulation, the  coordinatewise comultiplication). The reduced
semigroup bialgebra (i.e., the quotient by the one-dimensional ideal
generated by the zero of the semigroup)
is simply
the matrix bialgebra   $M_n(\Bbb C)$
(with the Kronecker comultiplication). As was noted above, both algebras
are nonunital involutive bialgebras. The semigroup
${\cal I}^1_n$ is a subsemigroup of the symmetric inverse semigroup
${\cal I}_n$ mentioned above.

\begin{theorem}
The semigroup algebra of a finite inverse unital semigroup
is a semisimple cocommutative involutive bialgebra. Analogously, the dual
semigroup algebra ${\Bbb C}(S)$ of a finite inverse unital semigroup $S$
is a commutative involutive bialgebra. Conversely, every finite-dimensional
semisimple cocommutative (in the dual version, commutative)
involutive bialgebra is isomorphic (as an involutive bialgebra)
to the semigroup algebra (respectively, the dual semigroup algebra)
of  a finite inverse unital semigroup.
\end{theorem}

This statement can be generalized
to inverse nonunital semigroups: one should only notice that the semigroup
algebra is not exactly a bialgebra, since either it has no counit or
the counit does not define a homomorphism into the field.

\begin{proof}
Let ${\Bbb C}(S)$ be the semigroup algebra over
$\Bbb C$ of a finite inverse semigroup $S$. Oganesyan's theorem
\cite{O} says that the semigroup algebra of every finite inverse semigroup
is semisimple. Hence we only have to check the conditions related to
the involution; namely, that
$(ab)^*=b^*a^*$. This is indeed true (see \cite[Lemma~1.18]{KP}).
The other properties of the involution and coinvolution are obvious.
If the semigroup $S$ has no two-sided unit, then there is no
corresponding homomorphism into $\Bbb C$, so that this is in fact
a weakened
involutive bialgebra.

To prove the converse, it is more convenient to consider the dual
formulation and the dual algebra $A'$. Since $A'$ is commutative and semisimple,
it follows that it is isomorphic, as an algebra with involution, to
${\Bbb C}^n$. Since the comultiplication
$\Delta$ is multiplicative, the image of an idempotent in
$A$ with respect to  $\Delta$ is an idempotent in $A\bigotimes A$,
hence the comultiplication determines a binary operation on the spectrum
of the algebra  ${\Bbb C}^n$, i.e., on
$\{1,2, \dots, n\}$. By the coassociativity of $\Delta$, this operation
defines on the spectrum the structure of a finite associative semigroup,
possibly without unit and possibly with zero. Since the involution commutes with
the comultiplication, it determines an involution on the semigroup.
This semigroup is inverse, as follows from the fact that every finite
subsemigroup of complex matrices closed under the ordinary involution of matrices
is inverse. In order to check this, one may use the following criterion
of being inverse: A semigroup with involution is inverse if and only if
the subsemigroup of idempotents is commutative
(see \cite[Theorem~1.17]{KP}); alternatively, one may use the analog of Cayley's theorem mentioned above
(see also \cite{Adv}).
The existence of a counit in the bialgebra is equivalent to the
existence of a unit in the semigroup.
\end{proof}

We have given a characterization of the class of semisimple cocommutative
finite-di\-men\-sion\-al involutive bialgebras. In contrast to group bialgebras,
semigroup bialgebras of inverse semigroups in general are not Hopf
algebras, since they have no antipode. But they have an ``almost antipode'':
given the semigroup algebra of an inverse semigroup, consider
the linear extension $S$ of the operation $a \rightarrow a^{-1}$ of taking the inverse;
this operator satisfies a condition that is
different from the ordinary condition on an antipode. Namely, the left and right
convolutions of $S$ with the identity map,
$S\star {\rm id}$ and ${\rm id} \star S$, are projections
to the commutative subalgebra
of idempotents rather than to the
one-dimensional subspace of scalars, as must be the case
for the antipode in a Hopf algebra. Thus the class of bialgebras we have described
is closest to the class of Hopf algebras.\footnote{Note that the theory
of representations of inverse semigroups by partially isometric
operators is also very close to the theory of unitary representations of
groups; besides, the semigroup bialgebra of a finite inverse semigroup
with unit is {\it Plancherel} in the sense of \cite{Ke1, VEP};
all this, in particular, allows one to simplify the proofs of
many facts, for example, the theorem on semisimplicity and other
properties of inverse semigroups mentioned above.}

\section{Positive 2-algebras}

\subsection{Definitions}

Theorem~1 motivates the following basic definition:

\begin{definition}
A semisimple 2-algebra $A$ over $\Bbb C$ with multiplication
 $\delta:A\bigotimes A\rightarrow A$, unit $\epsilon$,
involution $\sharp$, comultiplication $\Delta:A\rightarrow A\bigotimes
 A$, counit $\varepsilon$, and coinvolution $\flat$ is called
a \textbf{positive 2-algebra} if the operations are related as follows:

{\rm (1)} The multiplication and the comultiplication are positive, i.e., the
multiplication $\delta$
preserves the cone $K^{\flat}$, and the comultiplication
$\Delta$ preserves the cone
$K^{\sharp}$:  $\delta(K^{\flat}\bigotimes
K^{\flat})\subset K^{\flat}$, $\Delta K^{\sharp} \subset K^{\sharp}\bigotimes
K^{\sharp}$.

{\rm (2)} As in the definition of involutive bialgebras,
the involution $\sharp$ commutes with  the comultiplication,
and the coinvolution $\flat$ commutes with the multiplication:
$\Delta(x^{\sharp})=\{\Delta(x)\}^{\sharp \bigotimes \sharp}$,
$x^{\flat}\cdot y^{\flat}=(x\cdot y)^{\flat}$.
\end{definition}

The notion of a positive 2-algebra extends the notion of an involutive
bialgebra, replacing the condition that
 {\it the multiplication and the comultiplication should be multiplicative}
by the weaker condition that {\it both operations should be positive}.
In terms of duality of algebras, this notion was introduced in \cite{V}.
Theorem~1 implies the following corollary:

 \begin{corollary}
An involutive bialgebra satisfies conditions
{\rm (1), (2)} and hence is a positive 2-algebra.
 \end{corollary}

Indeed, condition~(1) follows from the multiplicativity of the multiplication and comultiplication
for bialgebras, and condition~(2) is contained in the definition of
involutive bialgebras.

Condition~(1) for positive 2-algebras cannot be replaced by the condition
that only one operation should be positive: here
the positivity of one operation does not, in general, imply the positivity of
the other one, as can be seen from simple examples
(see, e.g.,  \cite{BVEP,VEP}); this distinguishes positive 2-algebras
from bialgebras. Recognizing whether or not a finite-dimensional 2-algebra
is positive, given the finite-dimensional tensors of the
multiplication and comultiplication, is
an  $NP$-complete problem (in contrast to checking whether
a given 2-algebra is a bialgebra).

If we need  to emphasize the existence of all the structures mentioned above
in a bipositive algebra, we will write
 $$A=A(\delta,\epsilon,\sharp,\Delta,\varepsilon,\flat).$$

The most important class of involutive
bialgebras and positive 2-algebras is
obtained by adding another two dual
requirements relating both structures:

 \begin{definition}
A positive 2-algebra is called a
homogeneous positive 2-algebra if, in addition to conditions
{\rm (1), (2)}, it satisfies the condition

{\rm (3)} the counit is a positive homomorphism of the algebra into the field:
$\varepsilon (x\cdot y)=
 \varepsilon (x)\varepsilon (y)$; and the coproduct of the unit
is the unit in the tensor product:
 $\delta(\epsilon)=\epsilon\bigotimes\epsilon$. (These two conditions
 are symmetric with respect to the transition to dual algebras
 and are automatically satisfied for bialgebras.)
\end{definition}

The homogeneity condition (called so later in
\cite{VEP}) is exactly the condition of positive duality of pairs of algebras
stated in the original paper \cite{V}. It means that not only the cone
of positive (copositive) elements is closed under the comultiplication
(multiplication), but the compact sets of states (normalized positive
elements) are also closed under these operations.
This makes it possible to develop the so-called {\it geometric theory
of duality of packets},
or {\it block simplices},
i.e., pairs of compact affine semigroups of states
on algebras with involution, see
\cite{V, Ke1, Ke2}. {\it In what follows, unless otherwise stated,
we consider only homogeneous positive 2-algebras.\footnote{Without
going into details, we observe that, similarly to the fact
proved in \cite{VEP}, the category of homogeneous finite-dimensional semisimple
positive 2-algebras with the additional Plancherel condition (which is
a condition of general position) is equivalent to the category
of finite-dimensional algebras in Plancherel duality
in the sense of \cite{V,Ke2,VEP} and the category of positive
$C$-algebras \cite{Ban}.}} The recent paper \cite{Bu} (see the references
therein) contains
a description of the group algebras of
$n$-valued groups ($n$-Hopf algebras). This description is given in terms
of Frobenius homomorphisms. Since the group algebras of involutive $n$-valued
groups are a particular case of positive 2-algebras (namely, the case in which
the structure constants are rational), it would be interesting to
extend this description to all positive 2-algebras.

Somewhat freely, the difference between bialgebras and positive
homogeneous 2-algebras can be expressed as follows: the multiplication (and
comultiplication) in bialgebras is ``deterministic,''
while in positive
2-algebras it is ``probabilistic'' or multivalued
(see Sec.~3).

\subsection{Algebras in positive duality (\cite{V, K1})}

The notion of a positive 2-algebra is a paraphrase of the notion of
a pair of algebras in positive duality introduced in \cite{V}.
Let us give a very brief definition of this notion.
Let  $A$ and $B$ be two involutive unital algebras over
$\Bbb C$, and assume that there is a nondegenerate complex pairing
$<A,B>$ between these algebras. {\it The pair of algebras $A$ and $B$
is in a homogeneous positive duality} if the convex subset in $B$ of states
(= normalized positive definite functionals) of the algebra $A$ is stable
under the involution and multiplication in $B$, and, respectively, the convex
subset in $A$ of states of the algebra $B$ is stable under the involution
and multiplication in $A$; and, moreover, the unit of $A$ (respectively, $B$)
is the counit of $B$ (respectively, $A$). A positive duality is called
nonhomogeneous or weakened if only the cones of nonnegative elements
are stable in the same sense, but there is no condition on the units.
If we carry
the multiplication in $B$ (respectively, $A$) over to a comultiplication
in $A$ (respectively, $B$) in the standard way, then the obtained 2-algebra is a positive
2-algebra (respectively, the dual of a positive 2-algebra) in the sense of
the definition given in Sec.~5.1.

Sometimes, the language of duality is
more convenient than that of bialgebras.
For example, it allows one to develop
a meaningful study of the
finite-dimensional and infinite-dimensional geometry of pairs of {\it packets},
which was considered in
\cite{V} as a geometric theory of states on
$*$-algebras. A packet is a convex compactum of a certain kind
(a block simplex)
of states on a semisimple finite-dimensional
$*$-algebra endowed with a multiplication. The geometry of packets
is nontrivial even in the commutative case, where a packet is an affine
simplex with the structure of a semigroup with involution.
This geometry was much advanced in
\cite{Ke1,Ke2}. The corresponding definitions cover, as particular cases,
group algebras, algebras arising in algebraic combinatorics
(cellular algebras), and some new examples.
One may impose further restrictions
on the relation between the multiplication and the comultiplication
(the Plancherel property, etc.); on the author's initiative, they were
considered in detail in  \cite{Ke1, Ke2} and modified
in recent papers on algebraic combinatorics, multivalued groups, etc.,
see  \cite{BVEP, VEP, E, Bu}; here we do not consider these problems.

\section{Subobjects of involutive bialgebras and the statement
of the main problem}

\subsection{Strict subobjects}

Before formulating the central problem, we must define the notion
of a subobject in the category of involutive bialgebras. The author does not
know whether such a notion has been introduced for the category of bialgebras
or Hopf algebras. There are various possible definitions of a subobject.
We will consider two of them.
Let $A$ be  an involutive bialgebra,
and assume that $B$ is a unital subalgebra  of
$A$ (regarded as an algebra) closed under the involution
and coinvolution. If the condition
 $$\Delta B \subset B \bigotimes B
 \eqno{(*)}
 $$
is satisfied, then the comultiplication in $B$ is inherited from $A$, so that
it automatically
is associative and commutes with the involution.
In this case, $B$ is an involutive bialgebra
with the structures inherited from the bialgebra $A$, and it seems natural
to regard $B$ as a subobject in $A$. However, these conditions
are almost never satisfied, so that
the class of such subobjects
is too narrow and does not cover the most interesting and nontrivial applications.
Assume that instead of  ($*$), a much weaker condition is satisfied:

\begin{definition}
Let $A$ be a bialgebra and $B$ be
a unital subalgebra of $A$
closed under the involution
and coinvolution. It is called a strict subobject of the bialgebra $A$
or an admissible subalgebra of the involutive bialgebra $A$ if
there exists a coideal $J$ in $A$, also closed under the involution
and coinvolution, such that $A$ decomposes into the direct sum of
$B$ and $J$. The dual definition:
a strict subobject of a bialgebra is a subcoalgebra for which there exists
an ideal that decomposes the bialgebra into the direct sum of the
subcoalgebra and the ideal.
 \end{definition}

Under the natural identification of the quotient coalgebra
$A/J$ of the coalgebra $A$ by the coideal $J$ with the subalgebra
$B$, the latter obtains the structure of a coalgebra and turns into a
positive 2-algebra (see below); in general, it is not a bialgebra.
This definition can be stated in other words. Let $B$ be a unital subalgebra
of a bialgebra $A$ (regarded as an algebra)
closed under the involution
and coinvolution; assume that in $A$ there is a positive projection
$P=P_B$ from $A$ onto the subalgebra $B$\footnote{$P_B$ is
a positive expectation.}
 such that the map
 $$\Delta_{P_B}: B \rightarrow  (P\otimes P) \Delta, $$ i.e.,
the  $P$-projection of the comultiplication $\Delta$ to the subalgebra $B$,
determines a coassociative comultiplication on $B$, and assume that
the restriction of the counit of $A$ to $B$ is a counit with
respect to this comultiplication.

\begin{statement} For an admissible subalgebra $B$,
the multiplication $\delta|_B$ and the comultiplication
$\Delta_P$ introduced above are positive
with respect to the coinvolution and the involution, respectively.
\end{statement}

\begin{proof}
We need to check only the positivity of the comultiplication, but
a positive projection (= expectation)
sends positive operations (such as
a multiplicative comultiplication) to positive ones.
\end{proof}

This proposition, together with a direct check of condition~(3),
implies

\begin{corollary} An admissible subalgebra of an involutive bialgebra
is a homogeneous positive 2-algebra.
\end{corollary}

\medskip\noindent {\bf Remarks.}
1) The condition of the admissibility of a subalgebra
$B \subset A$ is very restrictive: although, as a rule,  for every subalgebra
there is a unique positive projection
$P_B$, the coassociativity of $(P_B\bigotimes P_B) \Delta$
holds only in special cases (see the next section).

2) The apparent asymmetry between the multiplication and the comultiplication
in our definition is dictated only by convenience considerations;
it is easy
to formulate an equivalent dual definition leading to the same class
of subobjects.
\medskip

Let us formulate the latter definition in terms of pairs of
algebras in positive duality, with the aim of giving a more extended definition
of a subobject in the same terms.

Consider an involutive bialgebra $\hat A$ over $\Bbb C$ as a pair of algebras
in duality $<A,A'>$; here $A$ is regarded as
$\hat A$ with the structure of an algebra, and
$A'$, the (algebraically) dual vector space, is endowed with the structure
of an algebra with the multiplication dual to the comultiplication in
$\hat A$. Let
$\hat B$ be a positive 2-algebra, and let
$<B,B'>$ be the pair of algebras in duality corresponding to
it in the same sense. By the above observations, in both cases
the duality of algebras is positive (in the second case, by definition).
The following proposition is a tautological reformulation of the
definition of a subobject given above.

\begin{statement}
Let $T$ be an isomorphic positive unit-preserving
embedding of the algebra $B$ into the algebra $A$
such that the restriction of the conjugate map
$T':A'\longrightarrow B'$ to some subalgebra $C$ of $A'$
is an isomorphism of $C$ and $B'$, and $T'$ is an expectation
to the subalgebra $C\subset A'$.
Then

{\rm 1)} under the identification of $C$ and $B'$
defined by $T'$, the multiplication in
$B'=C$ corresponds to the comultiplication in $B$;

{\rm 2)} $T\hat B$ is a subobject of $\hat A$.

Conversely, for every subobject $\hat B$ of the bialgebra $\hat A$,
the operator $T$ of embedding $B$ as an algebra into $A$
satisfies the above property.
\end{statement}

\subsection{Dynamical embeddings into involutive bialgebras:
nonstrict subobjects}

It turns out that the above definition of a subobject is too rigid
for the solution of the dilation problem; as we will see, the positive solution
is possible extremely rarely even in small dimensions, and the existence
conditions for such a solution are apparently difficult to formulate.
We will give a wider definition of a subobject, which could be called
{\it dynamical}. More precisely, we will define a new notion of
an embedding of a positive 2-algebra into an involutive bialgebra.
See \cite{Ca} for the definitions of a comodule and a coaction;
these definitions are easy to interpret in terms of pairs
of algebras in duality.

\begin{definition}
We say that a positive 2-algebra $\hat B$ can be nonstrictly embedded into a
bialgebra $\hat A$, or that its image is a nonstrict subobject of
$\hat A$, if there exists a positive isomorphic embedding of algebras
$T:\hat B \rightarrow \hat A$, preserving the unit, involution, and coinvolution,
such that the (left) coaction of $\hat B$ as a coalgebra on itself
can be extended from the algebra $TB$ to a (left) coaction of the whole
coalgebra  $\hat A$ on itself.
 \end{definition}

In the definition of a strict subobject, the comultiplication
in the subalgebra was being lifted (using a projection to the subalgebra)
to a comultiplication in the whole bialgebra; in the new definition, only the
coaction must be lifted, and there are no conditions on the way
in which this lifting should be realized. One may say that we embed
not a positive 2-algebra but rather this 2-algebra regarded as an algebra
with a coherent structure of a comodule.
The difference between the two notions of embedding
can be seen even for bicommutative algebras, where the multiplication
and the comultiplication are commutative. The dual statement can be obtained
by considering the dual coalgebras.

\subsection{Statement of the main problem}

Now we are ready to give an exact statement of the dilation problem, or the problem
of embedding positive 2-algebras into bialgebras.

\subsubsection{Strict version}
 {\it Let
 $$A=A(\delta,\epsilon,\sharp,\Delta,\varepsilon,\flat)$$
be a finite-dimensional homogeneous positive 2-algebra. When is there
a (perhaps, weakened) involutive bialgebra $\cal A $ such
that $A$ is isomorphic to some strict subobject $\cal B$ of $\cal A$?
If the answer is positive, then we will say that the corresponding involutive
bialgebra $\cal A$ is a strict \textbf{dilation} of the positive 2-algebra
$A$. When is there a minimal canonical strict dilation?}
\smallskip

Let us call this problem, as well as the next one,
the {\it problem of lifting a positive 2-algebra
into a bialgebra} or the {\it problem of strict dilation of a positive 2-algebra}
(in probabilistic terms, the {\it problem of derandomizing
a probabilistic comultiplication}). The homogeneity
allows us to interpret the dilation problem
as the problem of ``derandomizing''
a probabilistic comultiplication (see Sec.~3).

If we seek a solution of the problem
in the class of finite-dimensional involutive bialgebras,
then we certainly need to introduce additional restrictions on the
positive 2-algebra: it must have generators with
rational structure constants; the irrationality of the structure constants
in any generators
requires passing to infinite-dimensional algebras.

Taking into account Theorem~2, the main problem stated above
can be formulated in much more concrete terms:
{\it What finite-dimensional semisimple cocommutative homogeneous
positive 2-algebras with rational structure constants are isomorphic to
an admissible subalgebra of the dual semigroup algebra of an inverse
semigroup (in particular, a group)?
In other words, when there exists a strict dilation into the semigroup
algebra of an inverse semigroup of an arbitrary rational commutative
positive 2-algebra?} By analogy with groups, one could say that the
corresponding positive 2-algebras are of semigroup origin. The examples below
show that this extension of the class of algebras is not sufficient
for the positive solution of the problem, and for the two-dimensional
algebras, it gives nothing new compared with group bialgebras (see Sec.~6).

It presents no difficulty to state the problem in terms of
the structure constant tensor of the comultiplication in the
basis of minimal idempotents with respect to the multiplication, and then
formulate the problem itself as a problem of decomposing the
comultiplication tensors of a positive 2-algebra with respect to the Cayley
tensors of a group operation.
This formulation only emphasizes the difficulty of the problem, but
hardly provides a method of its solution.

\subsubsection{Nonstrict version}

Let us modify the statements of the previous section
by replacing the words ``strict subobject,'' ``strict dilation,'' etc.\
by ``nonstrict subobject,'' ``nonstrict dilation,'' etc. Then we will obtain
the nonstrict version of the problem. The only difference between the
two versions is that the latter uses
a wider notion of a subobject. This is a less rigid statement of the problem;
such a dilation will also be called {\it dynamical}. This notion is closer to
the notion of dilation in the sense of operator theory \cite{FN}.
However, the supply of positive 2-algebras for which a nonstrict dilation
is possible is much wider than in the first case; this can be seen already
in the case of bicommutative positive 2-algebras (see
the next section). The nonstrict version, in its geometric formulation,
is also illustrated by the example of Sec.~6.4.2.

\subsection{Examples of strict and nonstrict subobjects}

\subsubsection{Positive 2-algebras that are strict subobjects of
involutive bialgebras}

The complex group algebra ${\Bbb C}[G]$ of a finite group $G$,
regarded as the algebra of formal sums of group elements (or as the algebra
of complex measures on $G$ with the comultiplication that
carries such a measure from the group $G$ over to the diagonal of the
direct product $G\times G$ and the ordinary involution and coinvolution),
is an involutive bialgebra (and even a Hopf algebra with the antipode
$g \rightarrow g^{-1}$). The coproduct of a group element $g$
is the sum of group elements whose product
is equal to $g$.

The dual involutive bialgebra to the group algebra of a finite group
is a commutative involutive bialgebra: this is the space of all
complex functions on the group with pointwise multiplication,
diagonal comultiplication $\Delta:{\Bbb C}[G]\rightarrow {\Bbb
C}[G]\bigotimes{\Bbb C}[G]$, $(\Delta f)(g,h)=f(gh)$, involution
$(\sharp f)(g)=\bar f(g)$, and coinvolution
$(\flat f)(g)=\bar f(g^{-1})$.
The unit is
the function identically equal to one, and the counit is
the functional ``the value of a function at the identity
element of the group.''

A general example of an admissible
subalgebra (strict subobject) of a group bialgebra is as follows.

\begin{statement}
The subalgebra ${\Bbb C}[H\diagdown G\diagup H]\equiv {\cal
B}_H$ of double cosets of a subgroup
$H$ (sometimes called the Hecke algebra corresponding to the subgroup
$H$) in the group algebra ${\Bbb C}[G]$ is a positive 2-algebra that is
an admissible subalgebra of the bialgebra
${\Bbb C}[G]$ in the sense defined above.
\end{statement}
\begin{proof}
The projection $P$ from ${\Bbb C}[G]$ to ${\cal B}_H$
consists in averaging over the double
cosets of $H$; obviously, this projection is a positive expectation
(i.e. $P(xay)=xP(a)y$ for all $x,y \in {\cal B}_H$, $a \in {\Bbb C}[G]$) that
projects the comultiplication in ${\Bbb C}[G]$ to the natural
comultiplication
on the space of double cosets (dual to the multiplication of cosets). It is obvious that
the subalgebra ${\cal B}_H$ is closed under
the involution and coinvolution and contains the unit
and counit,
as well as that the required relations between the
operations are satisfied.
\end{proof}

The previous proposition remains literally true if we replace
the partition into the double cosets of a subgroup by an arbitrary {\it stable
partition} of the group, i.e., a partition such that
the space of functions constant on the partition blocks is an invariant
unital subalgebra of the group algebra; then, just as above,
this space is a positive 2-algebra with respect to the induced
structures. An example of a stable partition distinct from the partition
into the double cosets of a subgroup is the partition into the orbits of
some group of automorphisms  of the group.

This example can be extended to inverse semigroups. By Theorem~2,
the semigroup bialgebras of these semigroups
are involutive and semisimple, and one may consider
subobjects of these bialgebras in the above sense.
It is easy to check that the 2-algebra associated
with a stable partition of an inverse semigroup,
understood in the same sense as in the case of groups, is
an admissible subalgebra (subobject) of the semigroup bialgebra,
and the structure constants of such a 2-algebra are rational.
Conversely, the proof of Theorem~2 in fact contains the assertion that
an admissible subalgebra of the dual semigroup algebra is associated with
a stable partition of the inverse semigroup. Thus we obtain

\begin{corollary}
The class of strict subobjects of involutive finite-dimensional semisimple
bialgebras coincides with the class of positive 2-algebras generated by stable
partitions of finite inverse semigroups (in particular, groups).
\end{corollary}

Indeed, by Theorem~2, given such a bialgebra, we can first construct an
inverse semigroup and then a stable partition of this semigroup. In the
bicommutative case, the supplies of strict subobjects for group and semigroup
bialgebras coincide. Thus in this case the main problem (the dilation problem
or the problem of derandomizing a comultiplication) reduces to describing examples
of group origin (Schur algebras). It is these examples that were intensively
studied, but with another
interpretation and another terminology
(see \cite{VEP}),
in papers on the theory of association schemes.
However, this is only a small part of positive 2-algebras (or, in another
language, cellular algebras, $C$-algebras, etc.).
We emphasize that the theory of association schemes studies not merely
algebras of certain types, but their integer-valued representations.
The considerations of this work are algebraic; problems concerning
representations of algebras, especially integer-valued ones, should be
considered separately. Note that the strict dilation problem
can be stated in tensor terms (in terms of the
structure constants of
multiplication and comultiplication), similarly to the reformulation given below
of the nonstrict dilation problem for bicommutative 2-algebras in
terms of matrices of quasi-characters; this question will be considered
separately.

\subsubsection{Classification of two-dimensional algebras;
the difference between strict and nonstrict dilations}

The class of nonstrict subobjects is much wider than that
of strict ones; both have never been studied.
For want of space, we restrict ourselves
to an example, leaving details till another opportunity. Namely, we will
illustrate the difference between the two dilation problems with the example
of two-dimensional 2-algebras. Every two-dimensional positive
2-algebra is obviously bicommutative. Up to an isomorphism, it
can be described as follows. As a vector space, it is the two-dimensional
complex space ${\Bbb C}^2$. A structure of an algebra is defined by
generators, the unit $\textbf{1}$ of the algebra and an element
$\textbf{u}$, and the relation
$$\textbf{u}^2=(1-\lambda)\cdot\textbf{u}+\lambda \cdot \textbf{1}, \qquad \lambda \in [0,1].$$
The structure of a coalgebra is defined by the formulas
$$\Delta (\textbf{1})=\textbf{1}\otimes\textbf{1}, \qquad \Delta(\textbf{u}) = \textbf{u}\otimes \textbf{u}.$$
The involution and the coinvolution are the complex conjugation,
the counit $\epsilon$ is given by the formulas
$\varepsilon(\textbf{u})=-\lambda$,
$\varepsilon(\textbf{1})=1$.\footnote{In the dual description, the
generating idempotents are $\varepsilon$ and $v$,
and the comultiplication and the unit
are given by the formulas $\Delta \varepsilon=
\varepsilon \otimes
 \varepsilon$, $\Delta (v)= v \otimes 1 + 1\otimes v +(1-\lambda )v \otimes v$,
 $\textbf{1}(\varepsilon)=1$,  $\textbf{1}(v)=0$.} Denote the obtained
2-algebra by $A_{\lambda}$. It is easy to check that $A_{\lambda}$
is a positive 2-algebra, since all the axioms of Sec.~4 are satisfied. Note
that the multiplication determines the algebra of truncated polynomials
of second degree in the variable
$u$. In the standard basis of the space
${\Bbb C}^2$, which is regarded as a coalgebra, the unit takes the form
$\textbf{1}=(1,1)$ and the generator $\textbf{u}$ has the coordinates $(-\lambda,1)$.
It is easy to prove the following

\begin{statement}
Every two-dimensional positive semisimple 2-algebra is isomorphic to
one of the 2-algebras $A_{\lambda}$, $\lambda
\in [0,1]$. For $\lambda =1$, this is the group algebra of a group of second
order; for $\lambda=0$, the semigroup algebra of the semigroup
$\{<1,p> |p^2=p \}$. For $\lambda \in (0,1]$, the algebra $A_{\lambda}$
is Plancherel.
\end{statement}

\begin{theorem}
{\rm 1.} The positive 2-algebra $A_{\lambda}$ admits a strict dilation into an
involutive bialgebra (i.e., can be embedded as a strict subobject into
an involutive bialgebra, namely, into the group algebra of a finite
commutative group) if and only if
 $\lambda$ is a positive root of the equation
 $$z^2 -(2+\alpha)z + 1=0,$$
where $\alpha$ is of the form $\alpha=k\cdot \frac{(s-1)^2}{s}$
with $k,s$ arbitrary positive integers. Rational values of
$\lambda$ correspond to $k=1$
and are of the form $\lambda=n^{-1}$, $n\in \Bbb N$; for other $k$,
the values of $\lambda$ are irrational.

{\rm 2.} Assume that $\lambda$ is an arbitrary rational number; then
$A_{\lambda}$ admits a nonstrict dilation into the group algebra of
a finite commutative group.
\end{theorem}

For example, the algebra $A_{1/3}$ admits a nonstrict, but not a strict,
dilation into a bialgebra.

In the general case (see \cite{V, Ke1}), an $n$-dimensional homogeneous
bicommutative positive 2-algebra is determined by a matrix of
{\it quasi-characters}, i.e, a complex matrix satisfying
the following property: the coordinatewise product of any two rows
(columns) is a convex combination of rows (columns), and all the coordinates
of the first row and the first column are equal to one.
For the group bialgebra of a commutative group, this is the matrix
of characters; for the algebra
$A_{\lambda}$, this matrix has the form $\left(%
\begin{array}{cc}
  1 &       1    \\
  1 & -\lambda  \\
\end{array}%
\right)$, $\lambda \in [0,1]$.

The problem concerning nonstrict dilations of bicommutative positive 2-algebras
into bicommutative bialgebras reduces to the question whether
one can represent the matrix of quasi-characters as a coarse grain
of the matrix of characters of some commutative group
(or inverse semigroup).

\bigskip
Translated by N.~V.~Tsilevich.

\end{document}